\RequirePackage{fix-cm}
\documentclass[smallcondensed]{svjour3}     % onecolumn (ditto)
\smartqed  % flush right qed marks, e.g. at end of proof
\usepackage{graphicx}
\usepackage{amssymb}
\usepackage{amsmath,enumerate,bbm,color,verbatim,setspace}
\newtheorem{thm}{Theorem}
\newtheorem{lem}{Lemma}

\journalname{Statistical Inference for Stochastic Processes}
\begin{document}

\title{Optimal Filtering of Partially Observed Markov Processes with Gaussian Noise%\thanks{Grants or other notes
%about the article that should go on the front page should be
%placed here. General acknowledgments should be placed at the end of the article.}
}
%\subtitle{Do you have a subtitle?\\ If so, write it here}

%\titlerunning{Short form of title}        % if too long for running head

\author{Liubov A. Markovich  %etc.
}

%\authorrunning{Short form of author list} % if too long for running head

\institute{Liubov A. Markovich \at
              Institute of Control Sciences,
                Russian Academy of Sciences,
                Profsoyuznaya 65,
                117997 Moscow, Russia\\
              Tel.: 89168305103\\
              \email{kimo1@mail.ru}
}

\date{Received: date / Accepted: date}
% The correct dates will be entered by the editor

\maketitle

\begin{abstract}
The processing of stationary
random sequences under nonparametric uncertainty is given by a filtering problem
 when the signal distribution is unknown. A useful signal $(S_n)_{n\ge1}$ is assumed
  to be Markovian. This assumption allows
us to estimate the unknown $(S_n)$ using only an observable random sequence $(X_n)_{n\ge1}$. The equation of optimal filtering of such a
signal has been received by A.V. Dobrovidov. Our result states that
when the unobservable Markov sequence is defined by a linear equation with a
Gaussian noise, the equation of optimal filtering coincides both with the classical
Kalman filter and the conditional expectation defined by the Theorem on
normal correlation.
% \PACS{PACS code1 \and PACS code2 \and more}
% \subclass{MSC code1 \and MSC code2 \and more}
\end{abstract}

\section{Introduction}
\label{sec_0}
The problem of filtering of unknown signals from the mixture with noise has a wide range of applications including control of linear and nonlinear systems.
In the following we consider a partially observable Markov random sequence $(S_n,X_n)_{n\geq1}$, where a useful signal $S = (S_n)_{n\geq1}$  is unobservable and the sequence $X=(X_n)_{n\geq1}$ is observable.
The connection between these variables  is given by the following nonlinear (linear) expression
\begin{eqnarray}\label{24}
X_n&=&\varphi(S_n,\eta_n),
\end{eqnarray}
where $(\eta_n\in\mathbb{R})_{n\geq1}$ is an i.i.d random sequence, $(S_n)_{n\geq1}$ is a Markov sequence and $\varphi$ is some function. Realizations of random variables (r.v.s) $S_n\in\mathcal{S}_n\subseteq\mathbb{R}$ and $X_n\in\mathcal{X}_n\subseteq\mathbb{R}$ are denoted by $s_1^n=(s_1,\ldots,s_n)^T$ and $x_1^n=(x_1,\ldots,x_n)^T$, respectively.

\par In case when \eqref{24} has the recursive linear form
\begin{eqnarray}\label{0}S_n&=&aS_{n-1}+b\xi_n,\\
X_n&=&AS_n+B\eta_n\nonumber
\end{eqnarray}
where $S_n,X_n\in\mathbb{R}$ for all $n$; $\xi_n$ and $\eta_n$ are mutually independent r.v.s with the standard Gaussian distribution,
\begin{eqnarray*}   S_0\in\mathcal{N}(0,\widetilde{\sigma}^2),\quad \widetilde{\sigma}^2=\frac{b^2}{1-a^2}\\
S_n\in \mathcal{N}(0,1),\quad  n=1,2,3\ldots,
\end{eqnarray*}
coefficients $A,B,a,b$ are given real numbers and $|a|<1$, the Kalman filter is applied \cite{Kalman}.
However, the nonlinear models are more important for practice. The extended and unscented Kalman filter algorithms can be applied in this case, \cite{Crisan}, \cite{Julier}.

\par Another approach for nonlinear processes was proposed in \cite{Stratonovich:60}. With this respect, let us define the random sequence  $\vartheta_n=Q(S_n)$, where the r.v. $S_n$ is related to $\vartheta_n\in\Theta_n\subseteq\mathbb{R}$ by some one-to-one function $Q:\mathcal{S}_n\rightarrow\Theta_n$. The random sequence $(\vartheta_n)_{n\geq1}$ is also a Markov sequence.

\par To estimate $\vartheta_n$ the optimal Bayesian estimator in form of the conditional mean
\begin{eqnarray}\label{theta}\widehat{\vartheta}_n=\mathsf{E}(Q(S_n)|x_1^n)=\int\limits_{\mathcal{S}_n}Q(s_n)w_n(s_n|x_1^n)ds_n,
\end{eqnarray}
has been used. The  $w_n(s_n|x_1^n)$ is the posterior probability density function (pdf). It satisfies the Stratonovich's recurrence equation \cite{Stratonovich:66}
\begin{eqnarray}\nonumber w_{1}(s_{1}|x_1)&=&\frac{f(x_1|s_1)p(s_1)}{\int \limits_{\mathcal{S}_1}f(x_1|s_1)p(s_1)ds_1},\\
w_n(s_n|x_1^n)&=&\frac{f(x_n|s_n)}{f(x_n|x_1^{n-1})}\!\!\!\int\limits_{\mathcal{S}_{n-1}}p(s_n|s_{n-1})w_{n-1}(s_{n-1}|x_1^{n-1})ds_{n-1},\!\!\!\quad n\geq2.\label{1}
\end{eqnarray}
Here $p(s_n|s_{n-1})$ denotes the transition pdf of the Markov sequence $(S_n)_{n\geq1}$, $f(x_n|x_1^{n-1})$ and $f(x_n|s_n)$ denote conditional pdfs.
\par As the posterior pdf $w_n(s_n|x_1^n)$ depends on the
unknown prior distribution function $p(s_1)$ and the transition probability $p(s_n|s_{n-1})$ of the
Markov sequence $(S_n)_{n>1}$, we cannot use formula \eqref{theta} to estimate $\widehat{\vartheta}_n$.
To overcome this problem the optimal filtering equation (see Section \ref{sec_1}) was proposed in \cite{Dobrovidov:1983}.
\par We aim to prove the pairwise exact coincidences of the optimal filtering
equation in the form \eqref{th}, Kalman's filter and the conditional
expectation $\mathsf{E}(Q(S_n)|x_1^n)$ defined by the Theorem on normal correlation \cite{ShiryaevLiptser:2001}.  The latter
coincidences are shown to be valid when the unobservable Markov
sequence $(S_n)$ is defined by a linear equation with a Gaussian
noise.
Thus, the optimal filtering equation is nothing else but the Kalman filter in case of linear model \eqref{0}. However, for nonlinear processes
the optimal filtering equation provides the solution in contrast to the Kalman filter.
\par The Theorem on normal correlation requires the inverse covariance matrix. In case of the process \eqref{0} such matrix has a Toeplitz form \cite{Toeplitz}.
The explicit inversion of matrices from the Toeplitz class is considered in \cite{Trench:2001,Dow}.
In contrast to classical Toeplitz matrix we deal here with the modified matrix that differs from Toeplitz one by additional term on the diagonal. As an
auxiliary result we obtain the explicit inversion of such matrix. \\In \cite{ShiryaevLiptser:2001} a pseudo inverse matrix was used.
\par The paper is organized as follows. In Section \ref{sec_1} we remind the general equation of optimal filtering and its special case for the Gaussian pdf.
In Section \ref{sec_2} we obtain the conditional density $f(x_n|x_1^{n-1})$ and its derivative in explicit forms for the linear process \eqref{0} and the Gaussian pdf $f(x_n|s_n)$
(Theorem \ref{thm1}) and show their ratio to
estimate $E(S_n|x_1^n)$. In Section \ref{sec_3} the coincidence of the general filtration equation for the Gaussian pdf and the Kalman's
filter is derived (Theorem \ref{thm2}).
In Section \ref{sec_4} we find the explicit inverse covariance matrix  $D_{\vec{x}_n,\vec{x}_n}^{-1}$ and prove the coincidence  of the general filtration equation for the Gaussian pdf and the Theorem on normal correlation (Theorem \ref{thm3}).
All proofs are presented in the Appendices.
\section{Equations of optimal filtering}
\label{sec_1}
Motivated by the problem arising in Section \ref{sec_0}, we first transform \eqref{1} to a form which depends only on known variables.
\par Integrating \eqref{1} over $s_n$, we obtain
\begin{eqnarray*}\int\limits_{S_n}w_n(s_n|x_1^n)ds_n=\int\limits_{\mathcal{S}_n}\frac{f(x_n|s_n)}{f(x_n|x_1^{n-1})}\int\limits_{\mathcal{S}_{n-1}}p(s_n|s_{n-1})w_{n-1}(s_{n-1}|x_1^{n-1})ds_{n-1}ds_n.
\end{eqnarray*}
Furthermore, transferring $f(x_n|x_1^{n-1})$ to the left side of the latter equation we get
\begin{eqnarray}\label{2}f(x_n|x_1^{n-1})=\int\limits_{\mathcal{S}_n}f(x_n|s_n)\int\limits_{\mathcal{S}_{n-1}}p(s_n|s_{n-1})w_{n-1}(s_{n-1}|x_1^{n-1})ds_{n-1}ds_n.
\end{eqnarray}
Differentiation of \eqref{2} in $x_n$ leads to
\begin{eqnarray}\label{3}f'_{x_n}(x_n|x_1^{n-1})&=&\int\limits_{\mathcal{S}_n}f'_{x_n}(x_n|s_n)\int\limits_{\mathcal{S}_{n-1}}p(s_n|s_{n-1})w_{n-1}(s_{n-1}|x_1^{n-1})ds_{n-1}ds_n.
\end{eqnarray}
Let us further assume that the conditional pdf $f(x_n|s_n)$ belongs to the exponential family of distributions
\begin{eqnarray}\label{exponetfamily}f(x_n|s_n)=\widetilde{C}(s_n)h(x_n)\exp(T(x_n)Q(s_n)),
\end{eqnarray}
where $\widetilde{C}(s_n)$ is a normalization constant and $h(x_n),T(x_n),Q(s_n)$ are known functions. Its derivative in $x_n$ is given by
\begin{eqnarray*}f'_{x_n}(x_n|s_n)=f(x_n|s_n)\left(\frac{h'_{x_n}(x_n)}{h(x_n)}+T'_{x_n}(x_n)Q(s_n)\right).
\end{eqnarray*}
Substituting this into \eqref{3}, we can deduce that
\begin{eqnarray*}&&f'_{x_n}(x_n|x_1^{n-1})=\frac{h'_{x_n}(x_n)}{h(x_n)}f(x_n|x_1^{n-1})\\
&+&T'_{x_n}(x_n)\int\limits_{\mathcal{S}_n}f(x_n|s_n)Q(s_n)\int\limits_{\mathcal{S}_{n-1}}p(s_n|s_{n-1})w_{n-1}(s_{n-1}|x_1^{n-1})ds_{n-1}ds_n.
\end{eqnarray*}
Dividing the latter equation by $f(x_n|x_1^{n-1})$  and due to \eqref{1} we can write
\begin{eqnarray*}\frac{f'_{x_n}(x_n|x_1^{n-1})}{f(x_n|x_1^{n-1})}&=&\frac{h'_{x_n}(x_n)}{h(x_n)}+ T'_{x_n}(x_n)\int\limits_{\mathcal{S}_n}Q(s_n)w_{n}(s_{n}|x_1^{n})ds_n.
\end{eqnarray*}
Using \eqref{theta}, we can finally write that
\begin{eqnarray}\label{D}\mathsf{E}(Q(S_n)|x_1^n)\cdot T'_{x_n}(x_n)&=&%\frac{f'_{x_n}(x_n|x_1^{n-1})}{f(x_n|x_1^{n-1})}-\frac{h'_{x_n}(x_n)}{h(x_n)}=
\left(\ln\left(\frac{f(x_n|x_1^{n-1})}{h(x_n)}\right)\right)'_{x_n}.
\end{eqnarray}
This is a general filtration equation obtained in \cite{Dobrovidov:1983}.
Note that equation \eqref{D} does not contain the explicit probabilistic characteristics $p(s_1)$ and $p(s_n|s_{n-1})$ of the unknown sequence  $(S_n)$. This allows us to find the optimal estimator \eqref{theta} knowing only observable quantities of
$x_1^n$. Further, we shall call \eqref{D} as Dobrovidov's equation.
\par As an example of the exponential family \eqref{exponetfamily} we can take the Gaussian density
\begin{eqnarray}\label{f}f(x_n|s_n)=\frac{1}{\sqrt{2\pi}B}\exp\left(-\frac{(x_n-As_n)^2}{2B^2}\right).
\end{eqnarray}
Then the observation model is defined by the linear equation
\begin{eqnarray}\label{0_0}
X_n&=&AS_n+B\eta_n,
\end{eqnarray}
where $\eta_n$ is an i.i.d random sequence with Gaussian distribution and the coefficients $A$ and $B$ are given real numbers.
\par The pdf \eqref{f} relates to \eqref{exponetfamily}, where
\begin{eqnarray}\label{param1}\widetilde{C}(s_n)&=&\frac{1}{\sqrt{2\pi} B}\exp\left(-\frac{A^2s_n^2}{2B^2}\right),\quad h(x_n)=\exp\left(-\frac{x_n^2}{2B^2}\right),\\\nonumber
T(x_n)&=&x_n
,\quad
Q(s_n)=\frac{s_nA}{B^2}.
\end{eqnarray}
Substituting \eqref{param1} into \eqref{D}, we can write that
\begin{eqnarray}\label{th}\mathsf{E}(S_n|x_1^n)&=&\frac{B^2}{A}\frac{f'_{x_n}(x_n|x_1^{n-1})}{f(x_n|x_1^{n-1})}+\frac{x_n}{A}.
\end{eqnarray}
The latter formula is a special case of the general filtration equation  \eqref{D}. Furthermore, we need to obtain the conditional
density $f(x_n|x_1^{n-1})$ and its derivative.
\section{The conditional density $f(x_n|x_1^{n-1})$}\label{sec_2}
In this section we determine the conditional density \eqref{2} and its derivative in explicit forms. To this end, we consider
a partially observable  Markov sequence $(S_n,X_n)_{n\geq1}$ defined by the recursive linear equations \eqref{0}.
\par The following theorem holds.
\begin{thm}\label{thm1}The explicit form of the conditional density \eqref{2} is defined as
\begin{eqnarray}\label{frecur}&&f(x_n|x_1^{n-1})=\frac{1}{\sqrt{2\pi \sigma_n}}\exp\Biggl(-\frac{1}{2\sigma_n}\Bigl(x_n-A\mathcal{L}_{n-1}\Bigr)^2\Biggr),\quad \!\!n=2,3,\ldots,
\end{eqnarray}
where
\begin{eqnarray}\label{L}&&\mathcal{L}_n=\frac{Aa}{\sigma_{n-1}}\Biggl( x_{n-1}{\ae}_{n-1}+\frac{aB^2}{\sigma_{n-2}}\Bigl(x_{n-2} {\ae}_{n-2} + \frac{aB^2}{\sigma_{n-3}}\Bigl(x_{n-3}{\ae}_{n-3}+\ldots\nonumber\\
&+&\frac{aB^2}{\sigma_{2}}\Bigl(x_2{\ae}_2+x_1\frac{aB^2{\ae}_1}{\sigma_1}\underbrace{\Bigr)\ldots\Bigr)\Bigr)\Bigr)^2\Biggr)}_{n},\quad n=2,3,\ldots.
\end{eqnarray}
with
\begin{eqnarray}\label{param}{\ae}_1&=&\widetilde{\sigma}^2,\quad \sigma_1=B^2+A^2{\ae}_1,\\ \nonumber
{\ae}_n&=&\frac{B^2 a^2{\ae}_{n-1}+\sigma_{n-1} b^2}{\sigma_{n-1}},\quad
\sigma_n=B^2+A^2{\ae}_n, \quad n\geq2.
\end{eqnarray}
\end{thm}
The proof of Theorem \ref{thm1} is given in Appendix \ref{Ap1}.
\subsection{The ratio of the density derivative and the density}
Finally, we can find an explicit form of \eqref{th}.
To this end, we have to write the expression for the ratio of the derivative of the density and the density itself.
Using Theorem \ref{thm1}, it is straightforward to verify that
\begin{eqnarray}\label{ratio}&&\frac{f'_{x_n}(x_n|x_1^{n-1})}{f(x_n|x_1^{n-1})}=\frac{1}{\sigma_n}\Bigg(
\frac{A^2a}{\sigma_{n-1}}\Bigl(x_{n-1}{\ae}_{n-1}+\frac{aB^2}{\sigma_{n-2}}
\Bigl(x_{n-2} {\ae}_{n-2}+\\\nonumber
&+&\frac{aB^2}{\sigma_{n-3}}\Bigl(x_{n-3}{\ae}_{n-3}+\ldots\frac{aB^2}{\sigma_{2}}
\Bigl(x_2{\ae}_2+x_1\frac{aB^2{\ae}_1}{\sigma_1}\underbrace{\Bigr)\ldots\Bigr)\Bigr)\Bigr)}_{n}
-x_n\Bigg)
\end{eqnarray}
holds.
Substituting \eqref{ratio} into \eqref{th}, we can write
\begin{eqnarray}\label{12}\nonumber\mathsf{E}(S_n|x_1^n)
%=\\\nonumber
%&=&\frac{B^2}{A}\frac{1}{\sigma_n}\Big(\frac{A^2a}{\sigma_{n-1}}\Bigl(x_{n-1}{\ae}_{n-1}+\frac{aB^2}{\sigma_{n-2}}
%\Bigl(x_{n-2} {\ae}_{n-2}+\frac{aB^2}{\sigma_{n-3}}\Bigl(x_{n-3}{\ae}_{n-3}+\\\nonumber
%&+&\ldots\frac{aB^2}{\sigma_{2}}
%\Bigl(x_2{\ae}_2+x_1\frac{aB^2{\ae}_1}{\sigma_1}\underbrace{\Bigr)\ldots\Bigr)\Bigr)\Bigr)}_{n}
%-x_n\Big)+\frac{x_n}{A}\\\nonumber
&=&\frac{x_nA{\ae}_n}{\sigma_n}+\frac{x_{n-1}AaB^2{\ae}_{n-1}}{\sigma_{n-1}\sigma_n}+
\frac{x_{n-2}Aa^2B^4{\ae}_{n-2}}{\sigma_{n-2}\sigma_{n-1}\sigma_n}+\ldots+\\
&+&\frac{x_{2}Aa^{n-2}B^{2(n-2)}{\ae}_{2}}{\sigma_2\cdot\ldots\cdot\sigma_{n-3}\sigma_{n-2}\sigma_{n-1}\sigma_n}
+\frac{x_{1}Aa^{n-1}B^{2(n-1)}{\ae}_{1}}{\sigma_1\cdot\ldots\cdot\sigma_{n-3}\sigma_{n-2}\sigma_{n-1}\sigma_n}.
\end{eqnarray}
Using \eqref{L}, the ratio \eqref{ratio} can be represented by
\begin{eqnarray*}\frac{f'_{x_n}(x_n|x_1^{n-1})}{f(x_n|x_1^{n-1})}
&=&\frac{A\mathcal{L}_{n-1}-x_n}{\sigma_n}.
\end{eqnarray*}
Then Dobrovidov's equation \eqref{12} can be simplified to
\begin{eqnarray}\label{16}\mathsf{E}(S_n|x_1^n)&=&\frac{B^2}{A}\frac{A\mathcal{L}_{n-1}-x_n}{\sigma_n}+\frac{x_n}{A}
%=\\&=&\frac{x_n(1-B^2)}{A}+\frac{B^2S_{n-1}}{\sigma_n}
=\frac{Ax_n{\ae}_n}{\sigma_n}+\frac{B^2\mathcal{L}_{n-1}}{\sigma_n}.
\end{eqnarray}
Considering \eqref{16} one can represent \eqref{12} in a recursive form. We shall express $\mathsf{E}(S_{n+1}|x_1^{n+1})$ by $\mathsf{E}(S_n|x_1^n)$ using \eqref{16}.
As \eqref{L} can be represented as
\begin{eqnarray}\label{Srecur}&&\mathcal{L}_n=\frac{Aa}{\sigma_n}\Bigl(x_n{\ae}_n+\frac{\mathcal{L}_{n-1}B^2}{A}\Bigr).
\end{eqnarray}
it can be deduced that
\begin{eqnarray*}\mathsf{E}(S_{n+1}|x_1^{n+1})&=&\frac{Ax_{n+1}{\ae}_{n+1}}{\sigma_{n+1}}+\frac{B^2\mathcal{L}_{n}}{\sigma_{n+1}}=\\
&=&
\frac{Ax_{n+1}{\ae}_{n+1}}{\sigma_{n+1}}+\frac{B^2\left(\frac{Aa}{\sigma_n}\Bigl(x_n{\ae}_n+\frac{\mathcal{L}_{n-1}B^2}{A}\Bigr)\right)}{\sigma_{n+1}}\\
%&=&\frac{Ax_{n+1}{\ae}_{n+1}}{\sigma_{n+1}}+\frac{B^2Aax_n{\ae}_n}{\sigma_{n}\sigma_{n+1}}+\frac{B^4a\mathcal{L}_{n-1}}{\sigma_{n}\sigma_{n+1}}\\
&=&\frac{Ax_{n+1}{\ae}_{n+1}}{\sigma_{n+1}}+\frac{B^2Aax_n{\ae}_n}{\sigma_{n}\sigma_{n+1}}+\frac{B^2a}{\sigma_{n+1}}\left(\mathsf{E}(S_{n}|x_1^{n})-\frac{Ax_n{\ae}_n}{\sigma_n}\right).
\end{eqnarray*}
Therefore, Dobrovydov's equation \eqref{th} has a recursive form
\begin{eqnarray}\label{exp}\mathsf{E}(S_{n+1}|x_1^{n+1})&=&
\frac{Ax_{n+1}}{\sigma_{n+1}}{\ae}_{n+1}+\frac{B^2a}{\sigma_{n+1}}\mathsf{E}(S_{n}|x_1^{n}).
\end{eqnarray}
Later we shall use this form to prove
Theorem \ref{thm2}.
\section{The optimal filtering equation and its relation to Kalman filter}\label{sec_3}
Kalman filter for the linear
system \eqref{0} is defined by following recursive equations \cite{Dobrovidov:12}
\begin{eqnarray}\label{kalman}\mathsf{E}(S_{n+1}|x_1^{n+1})
%a\mathsf{E}(S_n|x_1^n)+\frac{Ab^2+a^2A\gamma_n}{B^2+A^2b^2+A^2a^2\gamma_n}(x_{n+1}-Aa\mathsf{E}(S_{n}|x_1^{n}))\\&=&
&=&\frac{Ab^2+a^2A\gamma_n}{B^2+A^2b^2+A^2a^2\gamma_n}x_{n+1}+\frac{aB^2\mathsf{E}(S_n|x_1^n)}{B^2+A^2b^2+A^2a^2\gamma_n},
\end{eqnarray}
\begin{eqnarray}\label{gamma}
\gamma_n&=&\frac{B^2(a^2\gamma_{n-1}+b^2)}{A^2(a^2\gamma_{n-1}+b^2)+B^2}
\end{eqnarray}
under the conditions
\begin{eqnarray*}\mathsf{E}(S_{1}|x_1)&=&\frac{A\widetilde{\sigma}^2}{A^2\widetilde{\sigma}^2+B^2}x_1,\quad
\gamma_1=\frac{B^2\widetilde{\sigma}^2}{A^2\widetilde{\sigma}^2+B^2}.\nonumber
\end{eqnarray*}
\par The following lemma holds.
\begin{lem}\label{lem0}
The parameters \eqref{gamma} are related to \eqref{param} by equation
\begin{eqnarray*}
\gamma_n&=&\frac{B^2{\ae}_n}{\sigma_n},
\end{eqnarray*}
where $B$ is given by \eqref{0}.
\end{lem}
\begin{thm}\label{thm2}When a partially observable Markov sequence $(S_n,X_n)_{n\geq1}$ is defined by  \eqref{0},
the equation of optimal filtering \eqref{th} is equivalent to the Kalman's filter \eqref{kalman}.
\end{thm}
Proofs of Lemma \ref{lem0} and Theorem \ref{thm2} are given in Appendices \ref{Ap2} and \ref{Ap3}.
\section{The Theorem on normal correlation}\label{sec_4}
In \cite{ShiryaevLiptser:2001} (Theorem 3.1, p.61) the Theorem on normal correlation has been obtained.
For the Gaussian vector $(\theta,\nu)$ the optimal estimate $\mathsf{E}(\theta|\nu)$
is defined by
\begin{eqnarray}\label{norm}\mathsf{E}(\theta|\nu)=\mathsf{E}(\theta)+D_{\theta\nu}D^{-1}_{\nu\nu}(\nu-\mathsf{E}(\nu)),\end{eqnarray}
where $\mathsf{E}(\theta)$ and
$\mathsf{E}(\nu)$ denote expectations and
\begin{eqnarray}\label{covmatr}D_{\theta\nu}&=&cov(\theta,\nu)=\|cov(\theta_i,\nu_j)\|,\quad 1\leq i\leq k, 1\leq j\leq l\\\nonumber
D_{\nu\nu}&=&cov(\nu,\nu)=\|cov(\nu_i,\nu_j)\|,\quad 1\leq i,j\leq l
\end{eqnarray}
are covariance matrices.
\par The Theorem on normal
 correlation \eqref{norm} contains the conditional mathematical expectation as the Dobrovydov's inequality \eqref{th}. It implies that \eqref{th} and \eqref{norm}
can be related. Therefore, we need to find how the covariance matrices \eqref{covmatr} can be expressed in terms of \eqref{0}.
\subsection{The covariance
matrices}
From \eqref{0} the following conditions
\begin{eqnarray*}\mathsf{E}(S_0)&=&0,\quad\mathsf{E}(S_0^2)=\frac{b^2}{1-a^2},\\
\mathsf{E}(\xi_n)&=&0,\quad\mathsf{E}(\eta_n)=0,\quad\mathsf{E}(X_n)=0,\quad n\geq1,\\
\mathsf{E}(\xi_n^2)&=&1,\quad\mathsf{E}(\eta^2_n)=1,\quad n\geq1
\end{eqnarray*}
follow.
Thus, using \eqref{0} we can write that $X_1=AS_1+B\eta_1$, and hence
$S_1=\frac{X_1-B\eta_1}{A}$ hold.
Then it follows
\begin{eqnarray*}S_2&=&aS_1+b\xi_2=\frac{a(X_1-B\eta_1)}{A}+b\xi_2,\\
S_3&=&aS_2+b\xi_3=\frac{a^2(X_1-B\eta_1)}{A}+ab\xi_2+b\xi_3.
%S_4&=&aS_3+b\xi_4=\frac{a^3(X_1-B\eta_1)}{A}+a^2b\xi_2+ab\xi_3+b\xi_4
\end{eqnarray*}
Let $S_n$ be defined as
\begin{eqnarray}\label{Sn}
S_n&=&\frac{a^{n-1}(X_1-B\eta_1)}{A}+a^{n-2}b\xi_2+\ldots+a^{n-(n-1)}b\xi_{n-1}+b\xi_{n}.
\end{eqnarray}
Then from \eqref{0}
\begin{eqnarray*}
S_{n+1}&=&aS_n+b\xi_{n+1}=\frac{a^{n}(X_1-B\eta_1)}{A}+a^{n-1}b\xi_2+\ldots+a^{n-n}b\xi_{n}+b\xi_{n+1}
\end{eqnarray*}
follows.
Thus, formula \eqref{Sn}  is true for any $n$ by induction.
\par Next, we can write down a similar formula for $X_n$, i.e. it holds
\begin{eqnarray*}
X_n&=&AS_n+B\eta_n\\
&=&a^{n-1}(X_1-B\eta_1)+a^{n-2}Ab\xi_2+\ldots+aAb\xi_{n-1}+Ab\xi_{n}+B\eta_n\\
&=&Aa^{n-1}(aS_0+b\xi_1)+a^{n-2}Ab\xi_2+\ldots+aAb\xi_{n-1}+Ab\xi_{n}+B\eta_n.
\end{eqnarray*}
Now we turn our attention to covariances.
\begin{lem}\label{lem01}
The following recursive formulas for the covariances
 \begin{eqnarray}\label{covxnxn}cov(X_n,X_n)&=&A^2{\ae}_1+B^2, \quad n\geq1,
\end{eqnarray}
 \begin{eqnarray}\label{covxmxn}cov(X_m,X_n)&=&A^2{\ae}_1a^{n-m},\quad n>m,\quad n\geq1,
\end{eqnarray}
\begin{eqnarray}\label{thetanxn}cov(S_n,X_n)&=&A{\ae}_1,\quad cov(S_n,X_m)=A{\ae}_1a^{n-m},\quad n>m,\quad n\geq1
\end{eqnarray}
hold.
\end{lem}
The proof of Lemma \ref{lem01} is given in Appendix \ref{Ap4}.
\par Next, combining \eqref{covxnxn}-\eqref{thetanxn}, we can finally write
the covariance matrices of any dimension $n$
\begin{eqnarray}\label{covmatrix}D_{\vec{X}_n,\vec{X}_n}
%&=&\left(
 %                                        \begin{array}{ccccc}
  %                                         cov(X_1,X_1) &cov(X_1,X_2) & cov(X_1,X_3)&\ldots& cov(X_1,X_n)\\
  %                                         cov(X_2,X_1) &  cov(X_2,X_2) & cov(X_2,X_3)&\ldots& cov(X_2,X_n) \\
   %                                        cov(X_3,X_1) & cov(X_3,X_2) & cov(X_3,X_3)&\ldots& cov(X_3,X_n)\\
    %                                        \vdots &  \vdots & \vdots& \ddots &\vdots\\
     %                                      cov(X_n,X_1) & cov(X_n,X_2) & cov(X_n,X_3)&\ldots& cov(X_n,X_n)\\
      %                                   \end{array}
       %                                \right)\\
                                       &=&\left(
                                         \begin{array}{ccccc}\nonumber
                                           var(X_1) &cov(X_1,X_2) & cov(X_1,X_3)&\ldots& cov(X_1,X_n)\\
                                           cov(X_2,X_1) &  var(X_2) & cov(X_2,X_3)&\ldots& cov(X_2,X_n) \\
                                           cov(X_3,X_1) & cov(X_3,X_2) & var(X_3)&\ldots& cov(X_3,X_n)\\
                                            \vdots &  \vdots & \vdots& \ddots &\vdots\\
                                           cov(X_n,X_1) & cov(X_n,X_2) & cov(X_n,X_3)&\ldots& var(X_n)\\
                                         \end{array}
                                       \right)\\
                                       &=&\left(
                                         \begin{array}{ccccc}
                                           A^2{\ae}_1+B^2 &A^2a{\ae}_1 & A^2a^2{\ae}_1&\ldots& A^2a^{n-1}{\ae}_1\\
                                           A^2a{\ae}_1 &  A^2{\ae}_1+B^2 & A^2a{\ae}_1&\ldots& A^2a^{n-2}{\ae}_1 \\
                                           A^2a^2{\ae}_1 &A^2a{\ae}_1 & A^2{\ae}_1+B^2&\ldots&A^2a^{n-3}{\ae}_1\\
                                            \vdots &  \vdots & \vdots& \ddots &\vdots\\
                                           A^2a^{n-1}{\ae}_1 & A^2a^{n-2}{\ae}_1 & A^2a^{n-3}{\ae}_1&\ldots& A^2{\ae}_1+B^2\\
                                         \end{array}
                                       \right),
%\\
                                      % &=&\left(
                                       %  \begin{array}{ccccc}
                                        %   \sigma_1 &A^2a{\ae}_1 & A^2a^2{\ae}_1&\ldots& A^2a^{n-1}{\ae}_1\\
                                         %  A^2a{\ae}_1 &  \sigma_1 & A^2a{\ae}_1&\ldots& A^2a^{n-2}{\ae}_1 \\
                                          % A^2a^2{\ae}_1 &A^2a{\ae}_1 & \sigma_1&\ldots&A^2a^{n-3}{\ae}_1\\
                                          %  \vdots &  \vdots & \vdots& \ddots &\vdots\\
                                          % A^2a^{n-1}{\ae}_1 & A^2a^{n-2}{\ae}_1 & A^2a^{n-3}{\ae}_1&\ldots& \sigma_1\\
                                        % \end{array}
                                      % \right),
\end{eqnarray}
\begin{eqnarray}\label{Dtheta}D_{S_n,\vec{X}_n}&=&
\left(
                                         \begin{array}{cccc}
                                           cov(S_n,X_1) &cov(S_n,X_2)&\ldots & cov(S_n,X_n) \\
                                         \end{array}
                                       \right)\\\nonumber
&=&
A{\ae}_1\left(
                                         \begin{array}{cccc}
                                           a^{n-1} &a^{n-2}&\ldots & 1 \\
                                         \end{array}
                                       \right).
\end{eqnarray}
Here $\vec{X}_n=X_1^n=(X_1,\ldots,X_n)^T$.
\par Matrix $D_{\vec{x}_n,\vec{x}_n}$ has to be inverted due to \eqref{norm}. It is not an easy problem to get an explicit matrix inversion.
Nevertheless, we further construct the inversion of our covariance matrix \eqref{covmatrix} for any dimension.
\subsection{The explicit inversion of Toeplitz matrix $D_{\vec{x}_n,\vec{x}_n}$}
The covariance matrix \eqref{covmatrix} is called a Toeplitz matrix \cite{Trench:2001} and can be represented as
\begin{eqnarray}\label{22}D_{\vec{X}_n,\vec{X}_n}&=&\left(
                                         \begin{array}{ccccc}
                                           c_1+B^2 &ac_1 & a^2c_1&\ldots& a^{n-1}c_1\\
                                           ac_1 &  c_1+B^2 & ac_1&\ldots& a^{n-2}c_1 \\
                                           a^2c_1 &ac_1 & c_1+B^2&\ldots&a^{n-3}c_1\\
                                            \vdots &  \vdots & \vdots& \ddots &\vdots\\
                                           a^{n-1}c_1 & a^{n-2}c_1 & a^{n-3}c_1&\ldots& c_1+B^2\\
                                         \end{array}
                                       \right)
\end{eqnarray}
where $c_1=A^2{\ae}_1$.
Hence, the covariance matrix \eqref{Dtheta} can be rewritten as
\begin{eqnarray}\label{23}D_{S_n,\vec{X}_n}&=&\frac{c_1}{A}\left(
                                         \begin{array}{ccccc}
                                           a^{n-1} &a^{n-2}&a^{n-3}&\ldots & 1 \\
                                         \end{array}
                                       \right).
\end{eqnarray}
We can represent \eqref{22} as follows
\begin{eqnarray*}D_{\vec{X}_n,\vec{X}_n}&=&(D_{\vec{X}_n,\vec{X}_n})_{B=0}+B^2\mathbf{I}
\end{eqnarray*}
where $\mathbf{I}$ is the identity matrix and $(D_{\vec{x}_n,\vec{x}_n})_{B=0}$ is determined by
\begin{eqnarray}\label{DB=0}(D_{\vec{X}_n,\vec{X}_n})_{B=0}&=&\left(
                                         \begin{array}{ccccc}
                                           c_1 &ac_1 & a^2c_1&\ldots& a^{n-1}c_1\\
                                           ac_1 &  c_1 & ac_1&\ldots& a^{n-2}c_1 \\
                                           a^2c_1 &ac_1 & c_1&\ldots&a^{n-3}c_1\\
                                            \vdots &  \vdots & \vdots& \ddots &\vdots\\
                                           a^{n-1}c_1 & a^{n-2}c_1 & a^{n-3}c_1&\ldots& c_1\\
                                         \end{array}
                                       \right).
\end{eqnarray}
Then the inverse matrix reads
\begin{eqnarray*}D_{\vec{X}_n,\vec{X}_n}^{-1}&=&((D_{\vec{X}_n,\vec{X}_n})_{B=0}+B^2\mathbf{I})^{-1}=\frac{1}{B^2}\left(\frac{1}{B^2}(D_{\vec{X}_n,\vec{X}_n})_{B=0}+\mathbf{I}\right)^{-1}.
\end{eqnarray*}
Using the formula $\left(\mathbf{P}+\mathbf{I}\right)^{-1} = \mathbf{P}^{-1} - \mathbf{P}^{-1}\left(\mathbf{I}+\mathbf{P}^{-1}\right)^{-1}\mathbf{P}^{-1}$, where $\mathbf{P}$ is a squared invertible matrix \cite{Gantmacher:1990}, we can write that
\begin{eqnarray}\label{inversed}&&D_{\vec{X}_n,\vec{X}_n}^{-1}=\frac{1}{B^2}\left(\frac{1}{B^2}(D_{\vec{X}_n,\vec{X}_n})_{B=0}+\mathbf{I}\right)^{-1}\\\nonumber
&=&(D_{\vec{X}_n,\vec{X}_n})_{B=0}^{-1}-B^2(D_{\vec{X}_n,\vec{X}_n})_{B=0}^{-1}
(\mathbf{I}+B^2(D_{\vec{X}_n,\vec{X}_n})_{B=0}^{-1})^{-1}(D_{\vec{X}_n,\vec{X}_n})_{B=0}^{-1}.
\end{eqnarray}
To find the inverse matrix $(D_{\vec{X}_n,\vec{X}_n})_{B=0}^{-1}$ one can use an algorithm from \cite{Trench:2001}.
\par Let $A_n$ be a squared, invertible  $n\times n$ matrix
 \begin{eqnarray*}A_n&=&\left(
                                         \begin{array}{ccccccc}
                                           1 &-a & 0 & 0&\ldots& 0&0\\
                                           0 & 1 & -a& 0&\ldots& 0 &0\\
                                           0 & 0 & 1& -a&\ldots&0&0\\
                                            \vdots &  \vdots & \vdots& \ddots&\ddots &\vdots&\vdots\\
                                           0 & 0 & 0&0& \ldots& 1&-a\\
                                           0 & 0 & 0& 0 &\ldots&0& 1\\
                                         \end{array}
                                       \right).
\end{eqnarray*}
Then it holds
\begin{eqnarray}\label{D1}(D_{\vec{X}_n,\vec{X}_n})_{B=0}\cdot A_n&=&\left(
                                         \begin{array}{ccccc}
                                           c_1 &0 & 0&\ldots& 0\\
                                           ac_1 &  \alpha_1 & ac_1&\ldots& 0 \\
                                           a^2c_1 &a\alpha_1 & \alpha_1&\ldots&0\\
                                            \vdots &  \vdots & \vdots& \ddots &\vdots\\
                                           a^{n-1}c_1 & a^{n-2}\alpha_1 & a^{n-3}\alpha_1&\ldots& \alpha_1\\
                                         \end{array}
                                       \right),
\end{eqnarray}
where $\alpha_1=c_1-a^2c_1$.
Therefore, by multiplying the left and the right sides of \eqref{D1} by the transposed matrix $A_n^T$ we can immediately write
\begin{eqnarray*}A_n^T(D_{\vec{X}_n,\vec{X}_n})_{B=0}A_n=diag(c_1,\alpha_1,\ldots,\alpha_1).\end{eqnarray*}
Multiplying the latter matrix from the left side by $(A^{-1}_n)^{T} $ and from the right-hand side by $A^{-1}_n$, we obtain
\begin{eqnarray*}(D_{\vec{X}_n,\vec{X}_n})_{B=0}=(A^{-1}_n)^{T}diag(c_1,\alpha_1,\ldots,\alpha_1)A^{-1}_n.\end{eqnarray*}
Hence, the inverse  matrix is given by
\begin{eqnarray*}(D^{-1}_{\vec{X}_n,\vec{X}_n})_{B=0}=A_ndiag(c_1^{-1},\alpha_1^{-1},\ldots,\alpha_1^{-1})A^{T}_n.\end{eqnarray*}
Finally, the inversion of the covariance matrix \eqref{DB=0} yields
\begin{eqnarray}\label{20}(D^{-1}_{\vec{X}_n,\vec{X}_n})_{B=0}&=&\left(
                                         \begin{array}{ccccc}
                                           c_1^{-1}+\alpha_1^{-1}a^2 &-\alpha_1^{-1}a & 0&\ldots& 0\\
                                           -\alpha_1^{-1}a  &  \alpha_1^{-1}+\alpha_1^{-1}a^2 & -\alpha_1^{-1}a&\ldots& 0 \\
                                          0 &-\alpha_1^{-1}a & \alpha_1^{-1}+\alpha_1^{-1}a^2&\ldots&0\\
                                            \vdots &  \vdots & \vdots& \ddots &\vdots\\
                                           0 & 0& 0&\ldots& \alpha_1^{-1}\\
                                         \end{array}
                                       \right)\nonumber\\
                                       &=&\frac{1}{c_1(1-a^2)}\left(
                                         \begin{array}{ccccc}
                                           1 &-a & 0&\ldots& 0\\
                                           -a  &1+a^2 &-a&\ldots& 0 \\
                                          0 &-a &1+a^2 &\ldots&0\\
                                            \vdots &  \vdots & \vdots& \ddots &\vdots\\
                                           0 & 0& 0&\ldots&  1\\
                                         \end{array}
                                       \right).
\end{eqnarray}
Next, using the notation $d_0-1=\frac{c_1(1-a^2)}{B^2}$ we can write
\begin{eqnarray*}\mathbf{I}+B^2(D^{-1}_{\vec{X}_n,\vec{X}_n})_{B=0}%&=&
 %\frac{1}{d_0-1}\left(
  %                                       \begin{array}{ccccc}
   %                                        d_0 &-a & 0&\ldots& 0\\
    %                                       -a  &d_0+a^2 &-a&\ldots& 0 \\
     %                                     0 &-a &d_0+a^2 &\ldots&0\\
      %                                      \vdots &  \vdots & \vdots& \ddots &\vdots\\
       %                                    0 & 0& 0&\ldots&  d_0\\
        %                                 \end{array}
         %                              \right)\\
                                       &=&\frac{a}{d_0-1}\left(
                                         \begin{array}{ccccc}
                                           \frac{d_0}{a} &-1 & 0&\ldots& 0\\
                                           -1  &\frac{d_0+a^2}{a} &-1&\ldots& 0 \\
                                          0 &-1 &\frac{d_0+a^2}{a} &\ldots&0\\
                                            \vdots &  \vdots & \vdots& \ddots &\vdots\\
                                           0 & 0& 0&\ldots&  \frac{d_0}{a}\\
                                         \end{array}
                                       \right).
\end{eqnarray*}
The latter matrix is a tridiagonal, symmetric matrix. In \eqref{inversed} we need its inverse. To this end, we
use the theory that was developed in \cite{Forseca:2007}, \cite{Usmani:1994}. Then we have
\begin{eqnarray}\label{ID^-1}(\mathbf{I}+B^2(D_{\vec{X}_n,\vec{X}_n})_{B=0}^{-1})^{-1}&=&\frac{d_0-1}{a}\left\{
                                                                                              \begin{array}{ll}
                                                                                                %(-1)^{i+j}(-1)^{j-1-i}\frac{\psi_{i-1}\varphi_{j+1}}{\psi_n}, & \hbox{if} \quad i\leq j \\
(-1)^{2j}\frac{\psi_{i-1}\varphi_{j+1}}{\psi_n}, & \hbox{if} \quad i\leq j \\
                                                                                                %(-1)^{i+j}(-1)^{i-1-j}\frac{\psi_{j-1}\varphi_{i+1}}{\psi_n}, & \hbox{if} \quad i>j.
(-1)^{2i}\frac{\psi_{j-1}\varphi_{i+1}}{\psi_n}, & \hbox{if} \quad i>j,
                                                                                              \end{array}
                                                                                            \right.
\end{eqnarray}
where $i,j=1,\ldots,n$ and $\psi,\varphi$ satisfy the following recurrence relations
\begin{eqnarray}\label{psi}\psi_m&=&\left(\frac{d_0+a^2}{a}\right)\psi_{m-1}-\psi_{m-2},\quad \mbox{for}\quad m=2,\ldots,n-1,
\end{eqnarray}
\begin{eqnarray}\label{psi2}
\psi_n&=&\frac{d_0}{a}\psi_{n-1}-\psi_{n-2},\quad \mbox{with initial conditions} \quad \psi_0=1, \psi_1=\frac{d_0}{a},\\\nonumber
\varphi_k&=&\left(\frac{d_0+a^2}{a}\right)\varphi_{k+1}-\varphi_{k+2},\quad \mbox{for}\quad k=n-1,\ldots,1\\\nonumber
&&\mbox{with initial conditions} \quad \varphi_{n+1}=1, \varphi_n=\frac{d_0}{a}.
\end{eqnarray}
\par Furthermore, $\psi_m=\varphi_{n+1-m}=\varphi_k,m=2,\ldots,n-1,\quad k=n-1,\ldots,1$. Thus, \eqref{ID^-1} can be expressed simply by
\begin{eqnarray}\label{21}&&(\mathbf{I}+B^2(D_{\vec{X}_n,\vec{X}_n})_{B=0}^{-1})^{-1}=\frac{d_0-1}{a\psi_n}\left\{
                                                                                              \begin{array}{ll}
                                                                                                %(-1)^{i+j}(-1)^{j-1-i}\frac{\psi_{i-1}\varphi_{j+1}}{\psi_n}, & \hbox{if} \quad i\leq j \\
\psi_{i-1}\psi_{n-j}, & \hbox{if} \quad i\leq j \\
                                                                                                %(-1)^{i+j}(-1)^{i-1-j}\frac{\psi_{j-1}\varphi_{i+1}}{\psi_n}, & \hbox{if} \quad i>j.
\psi_{j-1}\psi_{n-i}, & \hbox{if} \quad i>j.
                                                                                              \end{array}
                                                                                            \right.\nonumber\\
&=&\frac{d_0-1}{a\psi_n}\left(
     \begin{array}{cccccc}
       \psi_{n-1} & \psi_{n-2} & \psi_{n-3} & \ldots & \psi_{1} & 1\\
       \psi_{n-2} & \psi_{1}\psi_{n-2} & \psi_{1}\psi_{n-3} & \ldots & \psi_{1}^2 & \psi_{1} \\
       \psi_{n-3}& \psi_{1}\psi_{n-3} & \psi_{2}\psi_{n-3} & \ldots & \psi_{2}\psi_{1} & \psi_{2} \\
       \vdots & \vdots & \vdots & \ddots &\vdots & \vdots \\
       \psi_{1} & \psi_{1}^2 & \psi_{2}\psi_{1} & \ldots & \psi_{n-2}\psi_{1} & \psi_{n-2} \\
       1 & \psi_{1} & \psi_{2} & \ldots & \psi_{n-2} & \psi_{n-1}\\
     \end{array}
   \right).
\end{eqnarray}
%Then it is straightforward to verify that the inverse covariance matrix \eqref{inversed} is equal to
Replacing \eqref{20} and \eqref{21} into \eqref{inversed} one can obtain  the explicit inverse covariance matrix \eqref{covmatrix}.
\\As the product of \eqref{23} and \eqref{20} is given by
\begin{eqnarray*}&&D_{S_n,\vec{X}_n}\cdot(D_{\vec{X}_n,\vec{X}_n})_{B=0}=\frac{1}{A}\left(
                                         \begin{array}{ccccc}
                                          0 &0&0&\ldots & 1 \\
                                         \end{array}
                                       \right)
\end{eqnarray*}
then the product of the covariance matrices \eqref{23} and \eqref{inversed} is given by
\begin{eqnarray*}&&D_{S_n,\vec{X}_n}D_{\vec{X}_n,\vec{X}_n}^{-1}=\\
&=&-\frac{1}{Aa\psi_n}\Big(1-a\psi_1 \quad -a+(1+a^2)\psi_{1}-a\psi_{2} \quad -a\psi_{1}+(1+a^2)\psi_{2}-a\psi_{3}\quad \ldots \\
 &\ldots&-a\psi_{n-3}+(1+a^2)\psi_{n-2}-a\psi_{n-1} \quad -a\psi_{n-2}+\psi_{n-1}-a\psi_n\Big).
\end{eqnarray*}
Hence, the Theorem on normal correlation \eqref{norm} looks as follows
\begin{eqnarray}\label{13}&&\mathsf{E}(S_n|x_1^n)=D_{S_n,\vec{X}_n}D_{\vec{X}_n,\vec{X}_n}^{-1}\vec{x}_n=\\\nonumber
&=&\frac{a\psi_1-1}{Aa\psi_n}x_1-\frac{a\psi_{0}-(1+a^2)\psi_{1}+a\psi_{2}}{Aa\psi_n}x_2-\frac{a\psi_{1}-(1+a^2)\psi_{2}+a\psi_{3}}{Aa\psi_n}x_3 -\ldots \\\nonumber
&\ldots&-\frac{a\psi_{n-3}-(1+a^2)\psi_{n-2}+a\psi_{n-1}}{Aa\psi_n}x_{n-1}-\frac{a\psi_{n-2}-\psi_{n-1}+a\psi_n}{Aa\psi_n}x_{n}
\end{eqnarray}
\section{The Theorem on normal correlation and Dobrovidov's equation}\label{sec_5}
Formula \eqref{12} can be rewritten as follows
\begin{eqnarray}\label{mathE}&&\mathsf{E}(S_n|x_1^n)=\frac{Aa^{n-1}B^{2(n-1)}{\ae}_{1}}{\sigma_1\cdot\ldots\cdot\sigma_{n-2}\sigma_{n-1}\sigma_n}x_{1}+
\frac{Aa^{n-2}B^{2(n-2)}{\ae}_{2}\sigma_1}{\sigma_1\sigma_2\cdot\ldots\cdot\sigma_{n-2}\sigma_{n-1}\sigma_n}x_{2}
+\ldots\\\nonumber
&\ldots&+
\frac{AaB^2{\ae}_{n-1}\sigma_1\sigma_2\cdot\ldots\cdot\sigma_{n-2}}{\sigma_1\sigma_2\cdot\ldots\cdot\sigma_{n-3}\sigma_{n-2}\sigma_{n-1}\sigma_n}x_{n-1}+
\frac{A{\ae}_n\sigma_1\sigma_2\cdot\ldots\cdot\sigma_{n-2}\sigma_{n-1}}{\sigma_1\sigma_2\cdot\ldots\cdot\sigma_{n-3}\sigma_{n-2}\sigma_{n-1}\sigma_n}x_n
\end{eqnarray}
As we know from \eqref{psi}, $\psi_n,n=2,\ldots,N-1$ and $\psi_N$ are described by different formulas.
If $n=N$ is the number of the last element, we would mark the element as the last by $\widetilde{\psi}$. Then it holds
\begin{eqnarray*}
\widetilde{\psi}_N&=&\frac{d_0}{a}\psi_{N-1}-\psi_{N-2}.
\end{eqnarray*}
If the number of the last element is $n=N+1$, then we obtain
\begin{eqnarray*}
\psi_N&=&\frac{d_0+a^2}{a}\psi_{N-1}-\psi_{N-2}=\widetilde{\psi}_N+a\psi_{N-1}.
\end{eqnarray*}
Further, a similar representation can be written for the element $\psi_{N-1}$
\begin{eqnarray*}
\psi_{N-1}&=&\frac{d_0+a^2}{a}\psi_{N-2}-\psi_{N-3}=\widetilde{\psi}_{N-1}+a\psi_{N-2}
\end{eqnarray*}
and we get
\begin{eqnarray*}
\psi_N&=&\frac{d_0+a^2}{a}\psi_{N-1}-\psi_{N-2}=\widetilde{\psi}_N+a\widetilde{\psi}_{N-1}+a^2\psi_{N-2}.
\end{eqnarray*}
Repeating this procedure we obtain the following formulas
\begin{eqnarray}\label{14}\nonumber\psi_N&=&\sum\limits_{i=0}^{N-2}\widetilde{\psi}_{N-i}a^i+a^{N-1}\psi_{1},\\
\psi_{N-1}&=&\sum\limits_{i=1}^{N-2}\widetilde{\psi}_{N-i}a^{i-1}+a^{N-2}\psi_{1}
\end{eqnarray}
Then the last element $\widetilde{\psi}_{N+1}$ is the following
\begin{eqnarray}\label{psiN+1}\nonumber
\widetilde{\psi}_{N+1}&=&\frac{d_0}{a}\psi_{N}-\psi_{N-1}=\frac{d_0}{a}(\widetilde{\psi}_N+a\psi_{N-1})-\psi_{N-1}=
\frac{d_0}{a}\widetilde{\psi}_N+(d_0-1)\psi_{N-1}\\\nonumber
&=& \frac{d_0}{a}\widetilde{\psi}_N+(d_0-1)\sum\limits_{i=1}^{N-2}\widetilde{\psi}_{N-i}a^{i-1}+(d_0-1)a^{N-2}\psi_{1}\\
&=&\frac{d_0}{a}(\widetilde{\psi}_N+(d_0-1)a^{N-2})+(d_0-1)\sum\limits_{i=1}^{N-2}\widetilde{\psi}_{N-i}a^{i-1}.
\end{eqnarray}
The sum in the latter equation is not very convenient. Motivated by this problem we write
\begin{eqnarray*}
\widetilde{\psi}_N&=&\frac{d_0}{a}(\sum\limits_{i=1}^{N-2}\widetilde{\psi}_{N-i}a^{i-1}+a^{N-2}\frac{d_0}{a})
-\left(\sum\limits_{i=2}^{N-2}\widetilde{\psi}_{N-i}a^{i-2}+a^{N-3}\frac{d_0}{a}\right)\\
&=&\sum\limits_{i=1}^{N-2}\widetilde{\psi}_{N-i}a^{i-1}\left(\frac{d_0}{a}-\frac{1}{a}\right)+
\left(\frac{d_0}{a}\right)^2a^{N-2}-\frac{d_0}{a}a^{N-3}+\frac{\widetilde{\psi}_{N-1}}{a}
\end{eqnarray*}
where \eqref{14} was used. Hence, the sum reads
\begin{eqnarray}\label{sumN-1}&&\sum\limits_{i=1}^{N-2}\widetilde{\psi}_{N-i}a^{i-1}=
\frac{a}{d_0-1}\left(\widetilde{\psi}_{N}-\frac{\widetilde{\psi}_{N-1}}{a}-\left(\frac{d_0}{a}\right)^2a^{N-2}+\frac{d_0}{a}a^{N-3}\right)
\end{eqnarray}
Substituting \eqref{sumN-1} into \eqref{psiN+1}, we get
\begin{eqnarray}\label{psiN+12}\widetilde{\psi}_{N+1}&=&\frac{d_0}{a}\widetilde{\psi}_N+\frac{d_0}{a}(d_0-1)a^{N-2}\\\nonumber
&+&(d_0-1)\frac{a}{d_0-1}\left(\widetilde{\psi}_{N}-\frac{\widetilde{\psi}_{N-1}}{a}-\left(\frac{d_0}{a}\right)^2a^{N-2}+\frac{d_0}{a}a^{N-3}\right)\\\nonumber
&=&\widetilde{\psi}_N\left(\frac{d_0}{a}+a\right)-\widetilde{\psi}_{N-1}
\end{eqnarray}
\begin{lem}\label{lem1}If the last element $\psi_{n}$ has a number $n=N$, where $N\geq2$ is an integer number, then
\begin{eqnarray*}&&\widetilde{\psi}_N=\frac{(1-a^2)}{B^{2N}a^N}\prod\limits_{i=1}^{N}\sigma_i
\end{eqnarray*}
holds,
where $\widetilde{\psi}_{N}$ is the last element defined by \eqref{psi}.
\end{lem}
The proof of Lemma \ref{lem1} is given in  Appendix \ref{Ap5}.
\par Now we turn our attention to the numerators of \eqref{mathE}.
Let us introduce the following notations
\begin{eqnarray}\label{numer}C_{x_1}&=&Aa^{n-1}B^{2(n-1)}{\ae}_{1}=\frac{b^2}{1-a^2}Aa^{n-1}B^{2(n-1)},\\ \nonumber
C_{x_{i}}&=&Aa^{n-i}B^{2(n-i)}{\ae}_{i}\prod\limits_{j=1}^{i-1}\sigma_{j}, \quad i=2,\ldots,n.
\end{eqnarray}
Parameters \eqref{param} can be represented as
\begin{eqnarray}\label{kappa}{\ae}_n&=&\frac{B^2 a^2{\ae}_{n-1}+\sigma_{n-1} b^2}{\sigma_{n-1}}=\frac{B^2a^2}{A^2}+b^2-\frac{B^4a^2}{A^2\sigma_{n-1}},
\end{eqnarray}
\begin{eqnarray}\label{sigma}\sigma_n&=&B^2+A^2{\ae}_n%=B^2+A^2b^2+\frac{A^2B^2a^2({\ae}_{n-1}-B^2)}{A^2\sigma_{n-1}}\\\nonumber
=B^2+A^2b^2+B^2a^2-\frac{B^4a^2}{\sigma_{n-1}}, \quad n\geq2.
\end{eqnarray}
Hence, we can immediately write
\begin{eqnarray}\label{c_x}
 C_{x_{i}}&=&\frac{aB^2}{A}\left(B^2a^2+A^2b^2-\frac{B^4a^2}{\sigma_{i-1}}\right)\prod\limits_{j=1}^{i-1}\sigma_{j}\\\nonumber
&=&\frac{aB^2}{A}\left(1-\frac{B^2}{\sigma_{i}}\right)\prod\limits_{j=1}^{i}\sigma_{j},\quad i=2,\ldots,n.
\end{eqnarray}
Next, the following  lemmas can be proved.
\begin{lem}\label{lem2}The numerators of \eqref{13} can be represented as
\begin{eqnarray*}\frac{a\psi_{1}-1}{Aa}&=&\frac{Ab^2}{B^2a},\\
\frac{a\psi_{n-2}-\psi_{n-1}+a\psi_{n}}{Aa}
&=&\frac{Ab^2}{B^2a}\psi_{n-1},\quad n=2,\ldots,N-1,\\
\frac{a\psi_{N-2}-\psi_{N-1}+a\psi_{N}}{Aa}&=&\frac{Ab^2}{B^2a}\psi_{N-1}
\end{eqnarray*}
\end{lem}
The proof of Lemma \ref{lem2} is given in  Appendix \ref{Ap5}.
\begin{lem}\label{lem3}The numerators of \eqref{13} and \eqref{mathE} are related by
\begin{eqnarray}\label{cx1}\frac{a\psi_{1}-1}{Aa}&=&\frac{C_{x_{1}}(1-a^2)}{B^{2N}a^N},
\end{eqnarray}
\begin{eqnarray}\label{cx2}
\frac{a\psi_{n-2}-\psi_{n-1}+a\psi_{n}}{Aa}&=&\frac{C_{x_{n}}(1-a^2)}{B^{2N}a^N},\quad n=2,\ldots,N-1,
\end{eqnarray}
\begin{eqnarray}\label{cx3}
\frac{a\psi_{N-2}-\psi_{N-1}+a\psi_{N}}{Aa}&=&\frac{C_{x_{N}}(1-a^2)}{B^{2N}a^N},
\end{eqnarray}
where $C_{x_{i}}$ is defined by \eqref{numer}.
\end{lem}
\begin{thm}\label{thm3}The theorem on normal correlation \eqref{13} and Dobrovidov's equation \eqref{mathE} for the system \eqref{0} are coincided.
\end{thm}
The proofs of Lemma \ref{lem3} and Theorem \ref{thm3} are given in Appendices \ref{Ap7} and \ref{Ap8}.

\appendix
\section{Appendix section}\label{app}
\subsection{Proof of Theorem \ref{thm1}}\label{Ap1}
To prove \eqref{frecur} we  use mathematical induction. Thus, we have to prove that the statement of Theorem \ref{thm1} holds for $n=2$.
Using \eqref{2} we can write
\begin{eqnarray}\label{4}f(x_2|x_1)=\int\limits_{\mathcal{S}_2}f(x_2|s_2)\int\limits_{\mathcal{S}_{1}}p(s_2|s_{1})w_{1}(s_{1}|x_1)ds_{1}ds_2.
\end{eqnarray}
The conditional densities $f(x_1|s_1)$, $f(x_2|s_2)$ defined by \eqref{f} are Gaussian.
%\begin{eqnarray*}f(x_1|s_1)&=&\frac{1}{\sqrt{2\pi}B}\exp\left(-\frac{(x_1-As_1)^2}{2B^2}\right),\\ f(x_2|s_2)&=&\frac{1}{\sqrt{2\pi}B}\exp\left(-\frac{(x_2-As_2)^2}{2B^2}\right),\\
%p(s_1)&=&\frac{1}{\sqrt{2\pi}\widetilde{\sigma}}\exp\left(-\frac{s_1^2}{2\widetilde{\sigma}^2}\right),\quad\mbox{where}\quad \widetilde{\sigma}^2=\frac{b^2}{1-a^2}
%\end{eqnarray*}
Using the formula \eqref{1}, where
\begin{eqnarray*}
p(s_1)&=&\frac{1}{\sqrt{2\pi}\widetilde{\sigma}}\exp\left(-\frac{s_1^2}{2\widetilde{\sigma}^2}\right),\quad \widetilde{\sigma}^2=\frac{b^2}{1-a^2},
\end{eqnarray*}
we can write the posterior pdf as
\begin{eqnarray}\label{5}w_{1}(s_{1}|x_1)
%\frac{\frac{1}{\sqrt{2\pi}B}\exp\left(-\frac{(x_1-As_1)^2}{2B^2}\right)\frac{1}{\sqrt{2\pi}\widetilde{\sigma}}\exp\left(-\frac{s_1^2}{2\widetilde{\sigma}^2}\right)}{\int \limits_{\mathcal{S}_1} \left( \frac{1}{\sqrt{2\pi}B}\exp\left(-\frac{(x_1-As_1)^2}{2B^2}\right)\frac{1}{\sqrt{2\pi}\widetilde{\sigma}}\exp\left(-\frac{s_1^2}{2\widetilde{\sigma}^2}\right)\right)ds_1}\\\nonumber
&=&\frac{\exp\left(\frac{-(x_1-As_1)^2}{2B^2}-\frac{s_1^2}{2\widetilde{\sigma}^2}\right)}{\int \limits_{\mathcal{S}_1}\left( \exp\left(\frac{-(x_1-As_1)^2}{2B^2}-\frac{s_1^2}{2\widetilde{\sigma}^2}\right)\right)ds_1}%=\frac{\exp\left(\frac{-(x_1-As_1)^2}{2B^2}-\frac{s_1^2}{2\widetilde{\sigma}^2}\right)}{I_{den}}.
\end{eqnarray}
The integral in the  denominator of \eqref{5} can be reduced to the form
\begin{eqnarray*}I_{den}%&=&\int \limits_{\mathcal{S}_1}\left( \exp\left(\frac{-(x_1-As_1)^2}{2B^2}-\frac{s_1^2}{2\widetilde{\sigma}^2}\right)\right)ds_1\\
&=&
\int \limits_{\mathcal{S}_1} \left( \exp\left(-s_1^2\left(\frac{B^2+A^2\widetilde{\sigma}^2}{2b^2B^2}\right)+s_1\frac{Ax_1}{B^2}-\frac{x_1^2}{2B^2}\right)\right)ds_1.\\
\end{eqnarray*}
This is the Euler-Poisson integral that is known in the form
\begin{eqnarray}\label{euler}&&\int\limits_{-\infty}^{\infty}\exp(-x^2a^2+xb+c)dx=\frac{\sqrt{\pi}}{a}\exp\left(\frac{b^2}{4a^2}+c\right).
\end{eqnarray}
Thus, it is straightforward to verify that
\begin{eqnarray}\label{Iden}I_{den}&=&%\sqrt{\frac{2\pi B^2\widetilde{\sigma}^2}{B^2+A^2\widetilde{\sigma}^2}}\exp\left(-\frac{x_1^2}{2(B^2+A^2\widetilde{\sigma}^2)}\right)
\sqrt{\frac{2\pi}{\sigma}}\exp\left(-\frac{x_1^2}{2\sigma B^2\widetilde{\sigma}^2}\right),
\end{eqnarray}
where $\sigma=\frac{B^2+A^2\widetilde{\sigma}^2}{B^2\widetilde{\sigma}^2}$.
Substituting \eqref{Iden} into \eqref{5} we deduce the posterior pdf as
\begin{eqnarray}\label{w1}w_{1}(s_{1}|x_1)
&=&\sqrt{\frac{\sigma}{2\pi}}\exp\left(-\frac{\sigma}{2}\left(s_1-x_1\frac{A}{B^2\sigma}\right)^2\right).
\end{eqnarray}
Since the conditional density in the expression \eqref{4} is defined by
\begin{eqnarray*}p(s_2|s_{1})&=&\frac{1}{\sqrt{2\pi}b}\exp\left(-\frac{(s_2-as_1)^2}{2b^2}\right),
\end{eqnarray*}
 we can write using \eqref{w1} that
\begin{eqnarray}\label{6}&&\int\limits_{\mathcal{S}_{1}}p(s_2|s_{1})w_{1}(s_{1}|x_1)ds_{1}
%\\\nonumber
%&=&
%\int\limits_{\mathcal{S}_{1}}\frac{1}{\sqrt{2\pi}b}\exp\left(-\frac{(s_2-as_1)^2}{2b^2}\right)\!\!\!\sqrt{\frac{\sigma}{2\pi}}\exp\left(-\frac{\sigma}{2}\left(s_1-x_1\frac{A}{B^2\sigma}\right)^2\right)ds_{1}\\\nonumber
%&=&\sqrt{\frac{\sigma}{2\pi(\sigma b^2+a^2)}}\exp\left(\frac{-\left(Aax_1-\sigma B^2s_2\right)^2}{2B^4\sigma(\sigma b^2+a^2)}\right)\\\nonumber
%&=&\sqrt{\frac{1}{2\pi D}}\exp\left(-\frac{1}{2D}\left(x_1\frac{Aa}{\sigma B^2}-s_2\right)^2\right),
=
\frac{1}{\sqrt{2\pi {\ae}_2}}\exp\left(-\frac{\left(x_1\frac{Aa{\ae}_1}{\sigma_1}-s_2\right)^2}{2{\ae}_2}\right),
\end{eqnarray}
where the following notations are introduced%$D=\frac{\sigma b^2+a^2}{\sigma}$
\begin{eqnarray*}{\ae}_1&=&\widetilde{\sigma}^2,\quad \sigma_1=B^2+A^2{\ae}_1,\quad
{\ae}_2=\frac{B^2a^2{\ae}_1+\sigma_1b^2}{\sigma_1},\quad \sigma_2=B^2+A^2{\ae}_2.
\end{eqnarray*}
Using \eqref{6} in \eqref{4} and the  Euler-Poisson integral \eqref{euler} we deduce the conditional density for $n=2$ as
%\begin{eqnarray}\label{7}f(x_2|x_1)&=&\frac{1}{2\pi B\sqrt{D}}\int\limits_{\mathcal{S}_2}\left(-\frac{(x_2-As_2)^2}{2B^2}-\frac{1}{2D}\left(x_1\frac{Aa}{\sigma B^2}-s_2\right)^2\right)ds_2\\ \nonumber
%&=&\frac{1}{\sqrt{2\pi (A^2D+B^2)}}\exp\left(-\frac{\frac{1}{2DB^2}\left(x_2-\frac{x_1A^2a}{B^2\sigma}\right)^2}{\frac{A^2D+B^2}{DB^2}}\right)\\\nonumber
%&=&\frac{1}{\sqrt{2\pi\chi}}\exp\left(-\frac{1}{2\chi}\left(x_2-\frac{x_1A^2a}{B^2\sigma}\right)^2\right)\nonumber
%\end{eqnarray}
%where $\chi=A^2D+B^2$
\begin{eqnarray}\label{fn2}
f(x_2|x_1)&=&\frac{1}{\sqrt{2\pi\sigma_2}}\exp\left(-\frac{1}{2\sigma_2}\left(x_2-\frac{A^2a}{\sigma_1}x_1{\ae}_1\right)^2\right).
\end{eqnarray}
Thus, \eqref{6} and  \eqref{fn2} determine the basis of the mathematical induction.
\par The second step of the proof is to show that if the following formulas
\begin{eqnarray}\label{19}&&\int\limits_{\mathcal{S}_{n-1}}p(s_n|s_{n-1})w_{n-1}(s_{n-1}|x_1^{n-1})ds_{n-1}=\frac{1}{\sqrt{2\pi {\ae}_n}}\exp\Bigl(-\frac{(s_n-\mathcal{L}_{n-1})^2}{2{\ae}_n}\Bigr),\nonumber\\
&&f(x_n|x_1^{n-1})=\frac{1}{\sqrt{2\pi \sigma_n}}\exp\Bigl(-\frac{1}{2\sigma_n}\Bigl(x_n-A\mathcal{L}_{n-1}\Bigr)^2\Bigr)
\end{eqnarray}
for $n$ hold, where
\begin{eqnarray*}\mathcal{L}_{n-1}&=&\frac{Aa}{\sigma_{n-1}}\Bigl( x_{n-1}{\ae}_{n-1}+\frac{aB^2}{\sigma_{n-2}}\Bigl(x_{n-2} {\ae}_{n-2} +\\
&+&\frac{aB^2}{\sigma_{n-3}}\Bigl(x_{n-3}{\ae}_{n-3}+\ldots
\frac{aB^2}{\sigma_{2}}\Bigl(x_2{\ae}_2+x_1\frac{aB^2{\ae}_1}{\sigma_1}\underbrace{\Bigr)\ldots\Bigr)\Bigr)\Bigr)\Bigr)}_{n-1},
\end{eqnarray*}
where ${\ae}_n$ and $\sigma_{n}$ are defined by \eqref{kappa} and \eqref{sigma},
then also formulas \eqref{19} are valid for $n+1$.
\par For $n+1$ the posterior density is determined by
\begin{eqnarray*}&&w_n(s_n|x_1^n)=\frac{f(x_n|s_n)}{f(x_n|x_1^{n-1})}\int\limits_{\mathcal{S}_{n-1}}p(s_n|s_{n-1})w_{n-1}(s_{n-1}|x_1^{n-1})ds_{n-1}\\
&=&\frac{\frac{1}{\sqrt{2\pi}B}\exp\left(\frac{-(x_n-As_n)^2}{2B^2}\right)}{\frac{1}{\sqrt{2\pi \sigma_n}}\exp\Bigl(-\frac{1}{2\sigma_n}\Bigl(x_n-As_{n-1}\Bigr)^2\Bigr)}\frac{1}{\sqrt{2\pi {\ae}_n}}\exp\Bigl(-\frac{1}{2{\ae}_n}(s_n-\mathcal{L}_{n-1})^2)\Bigr)
\end{eqnarray*}
by its definition.
Thus, using \eqref{euler} and the latter formula we can rewrite \eqref{19} for the next step $n+1$, i.e.
\begin{eqnarray*}&&\int\limits_{\mathcal{S}_{n}}p(s_{n+1}|s_{n})w_{n}(s_{n}|x_1^{n})ds_{n}=\\
&=&
\int\limits_{\mathcal{S}_{n}}\frac{\frac{1}{\sqrt{8\pi^3{\ae}_n}bB}\exp\Big(-\frac{(s_{n+1}-as_{n})^2}{2b^2}-\frac{(x_n-A\mathcal{L}_n)^2}{2B^2}-\frac{1}{2{\ae}_n}\Big(s_n-\mathcal{L}_{n-1}\Big)^2\Big)}{\frac{1}{\sqrt{2\pi \sigma_n}}\exp\Bigl(-\frac{1}{2\sigma_n}\Bigl(x_n-As_{n-1}\Bigr)^2\Bigr)}ds_{n}ds_{n+1}\\
%&=&\sqrt{\frac{\sigma_n}{2\pi({\ae}_nB^2a^2+\sigma_nb^2)}}\exp\Bigg(\frac{\left(\frac{\mathcal{L}_{n-1}}{{\ae}_n}+\frac{Ax_n}{B^2}+\frac{as_{n+1}}{b^2}\right)^2}{2\left(\frac{\sigma_n}{B^2{\ae}_n}+\frac{a^2}{b^2}\right)}+\\
%&+&\left(\frac{(x_n-A\mathcal{L}_{n-1})^2}{2\sigma_n}-\frac{x_n^2}{2B^2}-\frac{\mathcal{L}_{n-1}^2}{2{\ae}_n}-\frac{s_{n+1}^2}{2b^2}\right)\Bigg)\\
&=&\frac{1}{\sqrt{2\pi{\ae}_{n+1}}}\exp\Bigg(-\frac{1}{2{\ae}_{n+1}}\Big(s_{n+1}-\Bigg(\frac{Aa}{\sigma_n}x_n{\ae}_n+\frac{B^2a\mathcal{L}_{n-1}}{\sigma_n}\Bigg)\Bigg)^2\Bigg)\\
%&=&\frac{1}{\sqrt{2\pi{\ae}_{n+1}}}\exp\Bigg(-\frac{1}{2{\ae}_{n+1}}\Big(s_{n+1}-\frac{Aa}{\sigma_n}\Bigg(x_n{\ae}_n+\frac{aB^2}{\sigma_{n-1}}
%\Bigl( x_{n-1}{\ae}_{n-1}+\\
%&+&\frac{aB^2}{\sigma_{n-2}}\Bigl(x_{n-2} {\ae}_{n-2} +\frac{aB^2}{\sigma_{n-3}}\Bigl(x_{n-3}{\ae}_{n-3}+\ldots
%\frac{aB^2}{\sigma_{2}}\Bigl(x_2{\ae}_2+x_1\frac{aB^2{\ae}_1}{\sigma_1}\underbrace{\Bigr)\!\ldots\!\Bigr)\!\!\Bigr)\!\!\Bigr)^2\!\Bigg)}_{n+1}\\
%&=&\frac{1}{\sqrt{2\pi{\ae}_{n+1}}}\exp\Bigg(-\frac{1}{2{\ae}_{n+1}}\Bigl(s_{n+1}-\frac{Aa}{\sigma_n}\Bigl(x_n{\ae}_n+\frac{\mathcal{L}_{n-1}B^2}{A}\Bigr)^2\Bigr)\\
&=&\frac{1}{\sqrt{2\pi{\ae}_{n+1}}}\exp\Bigg(-\frac{1}{2{\ae}_{n+1}}\Bigl(s_{n+1}-\mathcal{L}_{n}\Bigr)^2\Bigg)
\end{eqnarray*}
holds, where we use the notation \eqref{Srecur}.
Finally, we can write that
\begin{eqnarray*}&&f(x_{n+1}|x_1^{n})=\int\limits_{\mathcal{S}_{n+1}}f(x_{n+1}|s_{n+1})\int\limits_{\mathcal{S}_{n}}p(s_{n+1}|s_{n})w_{n}(s_{n}|x_1^{n})ds_{n}ds_{n+1}=\\
&=&\int\limits_{\mathcal{S}_{n+1}}\!\!\!\frac{\exp\left(\frac{-(x_{n+1}-As_{n+1})^2}{2B^2}-\frac{\left(s_{n+1}-\left(\frac{Aa}{\sigma_n}x_n{\ae}_n+\frac{B^2a\mathcal{L}_n}{\sigma_n}\right)\right)^2}{2{\ae}_{n+1}}
\right)}{\sqrt{4\pi^2{\ae}_{n+1}}B}
ds_{n+1}\\
&=&\frac{1}{\sqrt{2\pi\sigma_{n+1}}}\exp\left(-\frac{1}{\sigma_{n+1}}\left(x_{n+1}-A\left(\frac{Aa}{\sigma_n}x_n{\ae}_n+\frac{B^2a\mathcal{L}_n}{\sigma_n}\right)\right)^2\right)\\
%&=&\frac{1}{\sqrt{2\pi\sigma_{n+1}}}\exp\Big(-\frac{1}{\sigma_{n+1}}\Big(x_{n+1}-\frac{A^2a}{\sigma_n}(x_n{\ae}_n+\frac{aB^2}{\sigma_{n-1}}
%\Bigl( x_{n-1}{\ae}_{n-1}+\\
%&+&\frac{aB^2}{\sigma_{n-2}}\Bigl(x_{n-2} {\ae}_{n-2} +\frac{aB^2}{\sigma_{n-3}}\Bigl(x_{n-3}{\ae}_{n-3}+\ldots
%\frac{aB^2}{\sigma_{2}}\Bigl(x_2{\ae}_2+x_1\frac{aB^2{\ae}_1}{\sigma_1}\underbrace{\Bigr)\ldots\Bigr)\Bigr)^2\Bigr)}_{n+1}\\
&=&\frac{1}{\sqrt{2\pi\sigma_{n+1}}}\exp\Big(-\frac{1}{\sigma_{n+1}}\Big(x_{n+1}-A\mathcal{L}_n\Big)^2\Big).
\end{eqnarray*}
Since both the basis and the inductive step have been performed, by mathematical induction, the statement of Theorem \ref{thm1} holds for all integer $n>0$.
\subsection{Proof of Lemma \ref{lem0}}\label{Ap2}
As the basis we suppose that for $n=1$ the equation
\begin{eqnarray*}
\gamma_1&=&\frac{B^2\widetilde{\sigma}^2}{A^2\widetilde{\sigma}^2+B^2}=\frac{B^2{\ae}_1}{\sigma_1}
\end{eqnarray*}
is true.
\par  We have to show as the inductive step that if for $n$
\begin{eqnarray*}
\gamma_n&=&\frac{B^2{\ae}_n}{\sigma_n}
\end{eqnarray*}
holds, then it also holds for $n+1$.
\par By definition we get
\begin{eqnarray*}
\gamma_{n+1}&=&\frac{B^2(a^2\gamma_{n}+b^2)}{A^2(a^2\gamma_{n}+b^2)+B^2}=
\frac{B^2(a^2\frac{B^2{\ae}_n}{\sigma_n}+b^2)}{A^2(a^2\frac{B^2{\ae}_n}{\sigma_n}+b^2)+B^2}=\frac{B^2{\ae}_{n+1}}{\sigma_{n+1}}
\end{eqnarray*}
Since both the basis and the inductive step have been performed, by mathematical induction the statement of Lemma \ref{lem0} holds for all integer $n$.
\subsection{Proof of Theorem \ref{thm2}}\label{Ap3}
Note, that the denominator of the first term in \eqref{kalman} can be represented  as
\begin{eqnarray}\label{17}
B^2+A^2b^2+A^2a^2\gamma_n=B^2+A^2\left(\frac{B^2a^2{\ae}_n+\sigma_nb^2}{\sigma_n}\right)=\sigma_{n+1}
\end{eqnarray}
and its numerator as
\begin{eqnarray}\label{18}
Ab^2+a^2A\gamma_n=A\left(\frac{B^2a^2{\ae}_n+\sigma_nb^2}{\sigma_n}\right)=A{\ae}_{n+1}
\end{eqnarray}
Thus, \eqref{kalman} can be rewritten using \eqref{17}, \eqref{18} and Lemma \ref{lem0} as
\begin{eqnarray*}\mathsf{E}(S_{n+1}|x_1^{n+1})&=&
\frac{A{\ae}_{n+1}}{\sigma_{n+1}}x_{n+1}+\frac{B^2a}{\sigma_{n+1}}\mathsf{E}(S_{n}|x_1^{n})
\end{eqnarray*}
that coincides with \eqref{exp}. This implies that Dobrovidov's equation \eqref{th} under the condition \eqref{0} and the Kalman's filter are coincided.
\subsection{Proof of Lemma \ref{lem01}}\label{Ap4}
We have
\begin{eqnarray*}cov(X_n,X_n)&=&\mathsf{E}(X_n-\mathsf{E}(X_n))(X_n-\mathsf{E}(X_n))=\mathsf{E}(X_n^2)\\\nonumber
&=&\mathsf{E}(A^2a^{2(n-1)}(a^2S^2_0+b^2\xi^2_1)+a^{2(n-2)}A^2b^2\xi^2_2+\\\nonumber
&+&a^{2(n-3)}A^2b^2\xi^2_3+a^2A^2b^2\xi^2_{n-1}+A^2b^2\xi^2_{n}+B^2\eta^2_n)\\\nonumber
&=&A^2a^{2(n-1)}\left(a^2\frac{b^2}{1-a^2}+b^2\right)+a^{2(n-2)}A^2b^2+\\\nonumber
&+&a^{2(n-3)}A^2b^2+a^2A^2b^2+A^2b^2+B^2\\\nonumber
&=&A^2\!\!\left(a^{2(n-2)}b^2\left(\frac{a^2}{1-a^2}+1\right)+b^2a^{2(n-3)}+\ldots+a^2b^2+b^2\right)\!+B^2\\\nonumber
&=&\frac{A^2b^2}{1-a^2}+B^2=A^2{\ae}_1+B^2, \quad n=1,2,3\ldots
\end{eqnarray*}
Furthermore, it follows
 \begin{eqnarray*}cov(X_m,X_n)&=&\mathsf{E}(X_m-\mathsf{E}(X_m))(X_n-\mathsf{E}(X_n))=\mathsf{E}(X_n\cdot X_m)\\\nonumber
&=&A^2a^{n-1}a^{m-1}(a^2\mathsf{E}(S^2_0)+b^2\mathsf{E}(\xi^2_1))+a^{n-2}a^{m-2}A^2b^2\mathsf{E}(\xi^2_2)+\\\nonumber
&+&a^{n-3}a^{m-3}A^2b^2\mathsf{E}(\xi^2_3)+a^{m-m}a^{n-m}A^2b^2\mathsf{E}(\xi^2_{m})\\\nonumber
&=&b^2A^2\left(a^{n-2}a^{m-2}\left(\frac{a^2}{1-a^2}+1\right)+\ldots+a^{n-m}\right)\\\nonumber
&=&\frac{A^2b^2}{1-a^2}a^{n-m}=A^2{\ae}_1a^{n-m},\quad n>m,\quad n=1,2,3\ldots
\end{eqnarray*}
Similarly, we obtain the following two covariances
\begin{eqnarray*}cov(S_n,X_n)&=&\frac{Ab^2}{1-a^2}=A{\ae}_1,\\\nonumber
 cov(S_n,X_m)&=&\frac{Ab^2}{1-a^2}a^{n-m}=A{\ae}_1a^{n-m},\quad n>m,\quad n=1,2,3\ldots
\end{eqnarray*}
\subsection{Proof of Lemma \ref{lem1}}\label{Ap5}
Let us assume that for $\widetilde{\psi}_2$ the expression
\begin{eqnarray*}&&\widetilde{\psi}_2=\frac{(1-a^2)}{B^{4}a^2}\sigma_1\sigma_2
\end{eqnarray*}
is true.
\\We have to show as the inductive step that if for $\widetilde{\psi}_N$ the equation
\begin{eqnarray}\label{psiN}&&\widetilde{\psi}_N=\frac{(1-a^2)}{B^{2N}a^N}\prod\limits_{i=1}^{N}\sigma_i,
\end{eqnarray}
holds than the same holds for $\widetilde{\psi}_{N+1}$.
\par Let us substitute \eqref{psiN} into \eqref{psiN+12}. We get
\begin{eqnarray*}\widetilde{\psi}_{N+1}&=&\frac{(1-a^2)}{B^{2N}a^N}\left(\frac{d_0}{a}+a\right)\prod\limits_{i=1}^{N}\sigma_i
-\frac{(1-a^2)}{B^{2(N-1)}a^{N-1}}\prod\limits_{i=1}^{N-1}\sigma_i\\
&=&\frac{(1-a^2)\prod\limits_{i=1}^{N+1}\sigma_i}{B^{2(N+1)}a^{N+1}}\left(\frac{B^2a}{\sigma_{N+1}}\left(\frac{d_0}{a}+a\right)+\frac{B^4a^2}{\sigma_{N}\sigma_{N+1}}\right)\\
&=&\frac{(1-a^2)\prod\limits_{i=1}^{N+1}\sigma_i}{B^{2(N+1)}a^{N+1}}\left(\frac{\sigma_{N}\left(B^2a^2+A^2b^2+B^2-\frac{B^4a^2}{\sigma_{N}}\right)}{\sigma_{N}\sigma_{N+1}}\right).
\end{eqnarray*}
Finally, taking into account \eqref{sigma} we can write
\begin{eqnarray*}\widetilde{\psi}_{N+1}&=&\frac{(1-a^2)\prod\limits_{i=1}^{N+1}\sigma_i}{B^{2(N+1)}a^{N+1}}.
\end{eqnarray*}
\subsection{Proof of Lemma \ref{lem2}}\label{Ap6}
By definition we have $\psi_{1}=\frac{d_0}{a}=\frac{A^2b^2}{B^2a}$. Then it follows
\begin{eqnarray*}\frac{a\psi_{1}-1}{Aa}&=&\frac{Ab^2}{B^2a}.
\end{eqnarray*}
Using \eqref{psi} we can immediately write for the second numerator of \eqref{13}
\begin{eqnarray*}&&\frac{a\psi_{0}-(1+a^2)\psi_{1}+a\psi_{2}}{Aa}=
\frac{a\psi_{0}-(1+a^2)\psi_{1}+(d_0+a^2)\psi_{1}-a\psi_{0}}{Aa}\\
&=&\frac{-(1+a^2)\psi_{1}+\left(\frac{A^2b^2}{B^2}+1+a^2\right)\psi_{1}}{Aa}=\frac{Ab^2}{B^2a}\psi_{1}.
\end{eqnarray*}
Similarly, it can be done  for any $n=2,\ldots,N-1$. For example for $n=N-1$ we get
\begin{eqnarray}\label{xn-1}&&\frac{a\psi_{N-3}-(1+a^2)\psi_{N-2}+a\psi_{N-1}}{Aa}\\\nonumber
&=&
\frac{a\psi_{N-3}-(1+a^2)\psi_{N-2}+(d_0+a^2)\psi_{N-2}-a\psi_{N-3}}{Aa}\\\nonumber
&=&\frac{-(1+a^2)\psi_{N-2}+(\frac{A^2b^2}{B^2}+1+a^2)\psi_{N-2}}{Aa}=\frac{Ab^2}{B^2a}\psi_{N-2}.
\end{eqnarray}
For $n=N$ the numerator of \eqref{13} is different. Using \eqref{psi2} it can be deduced that
\begin{eqnarray*}&&\frac{a\psi_{N-2}-\psi_{N-1}+a\psi_{N}}{Aa}=
\frac{a\psi_{N-2}-\psi_{N-1}+d_0\psi_{N-1}-a\psi_{N-2}}{Aa}\\
&=&\frac{\psi_{N-1}(d_0-1)}{Aa}=\frac{Ab^2}{B^2a}\psi_{N-1}
\end{eqnarray*}
holds.
\subsection{Proof of Lemma \ref{lem3}}\label{Ap7}
For the first numerator of \eqref{13} we use \eqref{numer}. Then it is obvious that
\begin{eqnarray}C_{x_1}&=&\frac{Ab^2}{B^2a}\frac{a^NB^{2N}}{1-a^2}=\frac{a\psi_{1}-1}{Aa}\frac{a^NB^{2N}}{1-a^2}
\end{eqnarray}
holds.
Thus \eqref{cx1} follows.
For \eqref{cx2} it is enough to prove this statement for any $n=\{2,\ldots,N-1\}$. We shall show it for $n=N-1$.
As we know from Lemma \ref{lem2} the numerator of \eqref{13} for  $n=N-1$ is the following
\begin{eqnarray}\label{N-2}
\frac{a\psi_{N-3}-(1+a^2)\psi_{N-2}+a\psi_{N-1}}{Aa}&=&\frac{Ab^2}{B^2a}\psi_{N-2}.
\end{eqnarray}
Therefore, if we use the same technique as in Lemma \ref{lem1} we can represent $\psi_{N-2}$ as
\begin{eqnarray}\label{psiN-2}\psi_{N-2}&=&\widetilde{\psi}_{N-2}+a\psi_{N-3}=\sum\limits_{i=2}^{N-2}\widetilde{\psi}_{N-i}a^{i-2}+a^{N-3}\widetilde{\psi}_1\\\nonumber
&=&\sum\limits_{i=2}^{N-2}\widetilde{\psi}_{N-i}a^{i-2}+a^{N-4}d_0
\end{eqnarray}
Using the similar technique as in \eqref{psiN+12} we get
\begin{eqnarray*}\widetilde{\psi}_{N-1}&=&\frac{d_0}{a}\left(\sum\limits_{i=2}^{N-2}\widetilde{\psi}_{N-i}a^{i-2}+a^{N-3}\frac{d_0}{a}\right)-
\left(\sum\limits_{i=3}^{N-2}\widetilde{\psi}_{N-i}a^{i-3}+a^{N-4}\frac{d_0}{a}\right)\\
&=&\sum\limits_{i=2}^{N-2}\widetilde{\psi}_{N-i}a^{i-2}\left(\frac{d_0}{a}-\frac{1}{a}\right)+\left(\frac{d_0}{a}\right)^2a^{N-3}-\frac{d_0}{a}a^{N-4}+
\frac{\widetilde{\psi}_{N-2}}{a}.
\end{eqnarray*}
Next, expressing from the latter equation the sum
\begin{eqnarray*}&&\sum\limits_{i=2}^{N-2}\widetilde{\psi}_{N-i}a^{i-2}=\frac{a}{d_0-1}\left(\widetilde{\psi}_{N-1}-
\frac{\widetilde{\psi}_{N-2}}{a}-\left(\frac{d_0}{a}\right)^2a^{N-3}+\frac{d_0}{a}a^{N-4}\right)
\end{eqnarray*}
and  substituting it into \eqref{psiN-2} we can write
\begin{eqnarray*}\psi_{N-2}&=&\frac{a}{d_0-1}\left(\widetilde{\psi}_{N-1}-
\frac{\widetilde{\psi}_{N-2}}{a}-\left(\frac{d_0}{a}\right)^2a^{N-3}+\frac{d_0}{a}a^{N-4}\right)+d_0a^{N-4}\\
&=&\frac{a\widetilde{\psi}_{N-1}-\widetilde{\psi}_{N-2}}{d_0-1}=\frac{B^2(a\widetilde{\psi}_{N-1}-\widetilde{\psi}_{N-2})}{A^2b^2}\\
&=&\frac{B^2(1-a^2)}{A^2b^2}\left(\frac{a\prod\limits_{j=1}^{N-1}\sigma_{j}}{B^{2(N-1)}a^{N-1}}-\frac{\prod\limits_{j=1}^{N-2}\sigma_{j}}{B^{2(N-2)}a^{N-2}}\right)
\end{eqnarray*}
We substitute the latter results into \eqref{N-2}. Finally, we deduce that
\begin{eqnarray*}\frac{a\psi_{N-3}-(1+a^2)\psi_{N-2}+a\psi_{N-1}}{Aa}&=&\frac{(1-a^2)}{AB^{2(N-1)}a^{N-1}}\left(1-\frac{B^2}{\sigma_{N-1}}\right)\prod\limits_{j=1}^{N-1}\sigma_{j}\\
&=&C_{x_{N-1}}\frac{(1-a^2)}{B^{2N}a^{N}}
\end{eqnarray*}
holds, where  the definition \eqref{c_x} of $C_{x_{N-1}}$  is used.
Thus, the statement \eqref{cx2} is proved.
For the last case when $n=N$ the same results are valid. From Lemma \ref{lem2}  we get that
\begin{eqnarray}\label{psiN2}
\frac{a\psi_{N-2}-\psi_{N-1}+a\psi_{N}}{Aa}&=&\frac{Ab^2}{B^2a}\psi_{N-1}
\end{eqnarray}
holds.
From \eqref{14}, \eqref{sumN-1} it follows
\begin{eqnarray*}
\psi_{N-1}&=&\frac{a\widetilde{\psi}_{N}-\widetilde{\psi}_{N-1}}{d_0-1}=
\frac{B^2(1-a^2)}{A^2b^2}\left(\frac{a\prod\limits_{j=1}^{N}\sigma_{j}}{B^{2N}a^{N}}-\frac{\prod\limits_{j=1}^{N-1}\sigma_{j}}{B^{2(N-1)}a^{N-1}}\right)
\end{eqnarray*}
Substituting it into \eqref{psiN2} and using \eqref{c_x} we deduce
\begin{eqnarray*}
\frac{a\psi_{N-2}-\psi_{N-1}+a\psi_{N}}{Aa}&=&\frac{(1-a^2)}{AB^{2N}a^N}\left(1-\frac{B^2}{\sigma_N}\right)\prod\limits_{j=1}^{N}\sigma_{j}
=C_{x_{N}}\frac{(1-a^2)}{B^{2N}a^{N}}.
\end{eqnarray*}
\subsection{Proof of Theorem \ref{thm3}}\label{Ap8}
The assertion of the theorem  follows immediately from Lemmas \ref{lem1}  and \ref{lem2}.
\begin{acknowledgements}
L. A. M. acknowledges the financial support provided within the Russian Foundation for Basic Research, grant 13-08-00744. A.
L. A. M. would like to thank Prof. A.V. Dobrovidov and Prof. R.S. Liptser for reading a draft of this paper and giving useful comments.
\end{acknowledgements}

%
%If you'd like to thank anyone, place your comments here
%and remove the percent signs.
%

% BibTeX users please use one of
%\bibliographystyle{spbasic}      % basic style, author-year citations
%\bibliographystyle{spmpsci}      % mathematics and physical sciences
%\bibliographystyle{spphys}       % APS-like style for physics
%\bibliography{}   % name your BibTeX data base
\bibliographystyle{spbasic}
\bibliography{reference}

% Non-BibTeX users please use
%\begin{thebibliography}{}
%
% and use \bibitem to create references. Consult the Instructions
% for authors for reference list style.
%
%\bibitem{RefJ}
% Format for Journal Reference
%Author, Article title, Journal, Volume, page numbers (year)
% Format for books
%\bibitem{RefB}
%Author, Book title, page numbers. Publisher, place (year)
% etc
%\end{thebibliography}

\end{document}